\documentclass[11pt]{article}
\usepackage{minitoc}
\usepackage{amssymb}
\usepackage{amsmath}
\usepackage{mathrsfs}
\usepackage{amsfonts}
\usepackage{amsfonts,psfrag,graphicx}
\usepackage[a4paper]{geometry}
\usepackage{color}

\mtcselectlanguage{french} 
\newtheorem{defi}{Def\hskip 1pt inition }[section]
\newtheorem{lem}[defi]{Lemma}
\newtheorem{theo}[defi]{Theorem}
\newtheorem{prop}[defi]{Proposition}

\newtheorem{techlem}[defi]{Technical Lemma}
\newtheorem{rem}[defi]{Remark}
\newtheorem{assum}[defi]{Assumption}

\def\pn{\par\noindent}

\def\videbox{\mathbin{\vbox{\hrule\hbox{\vrule height1ex \kern.5em\vrule
height1ex}\hrule}}}
\def\1{1\!{\rm l}}

 \newcommand{\bigO}[1]{\ensuremath{\mathop{}\mathopen{}O\mathopen{}\left(#1\right)}}
 \newcommand{\smallO}[1]{\ensuremath{\mathop{}\mathopen{}o\mathopen{}\left(#1\right)}}
 
 \newcommand{\di}{\displaystyle} 
 \newcommand{\footnoteremember}[2]{
\footnote{#2}
\newcounter{#1}
\setcounter{#1}{\value{footnote}}
}
\newcommand{\footnoterecall}[1]{
\footnotemark[\value{#1}]
}

\begin{document}
\title{Sharp Large Deviations for empirical correlation coefficients}
\author{T.K.T. Truong \footnoteremember{IDP}{Institut Denis Poisson, Universit\'e d'Orl\'eans, Universit\'e de Tours, CNRS, Route de Chartres, B.P. 6759, 45067, Orl\'eans cedex 2, France}, M. Zani \footnoterecall{IDP}\footnote{Corresponding author: marguerite.zani@univ-orleans.fr}
}
\date{}
\maketitle
\pn

\begin{abstract}
In this paper, we study Sharp Large Deviations for empirical Pearson coefficients, i.e.  $r_n=\sum_{i=1}^n(X_i-\bar X_n)(Y_i-\bar Y_n)/\sqrt{\sum_{i=1}(X_i-\bar X_n)^2 \sum_{i=1}(Y_i-\bar Y_n)^2}$ or
 $\tilde r_n=\sum_{i=1}^n(X_i-\mathbb E(X))(Y_i-\mathbb E(Y))/\sqrt{\sum_{i=1}(X_i-\mathbb E(X))^2 \sum_{i=1}(Y_i-\mathbb E(Y))^2}\, $ (when the expectations are known).
Our framework is for random samples $(X_i,Y_i)$ either Spherical or Gaussian. In each case, we follow the scheme of Bercu et al. We also compute the Bahadur exact slope in the Gaussian case.
\end{abstract}

\maketitle
\section{Introduction}
The Beauvais--Pearson linear correlation coefficient between two real random variables $X$ and $Y$ is defined by
 \begin{align*}
\rho =\dfrac{\text{Cov}(X,Y)}{	\sqrt{\text{Var}(X) }	\sqrt{\text{Var}(Y)}},
\end{align*}
whenever this quantity exists. Such quantities were formally defined more than a century ago by Pearson \cite{PearsonHist, Pearson1}. The correlation describes the linear relation between two random variables. It is clear from Cauchy--Schwartz inequality that the absolute value of $\rho$ is less than or equal to $1$. Moreover, $\rho= \pm 1$ if and only if $X$ and $Y$ are linearly related. When $\rho= 0$  we say that $X$ and $Y$ are uncorrelated, i.e. linearly independent. 
The empirical counterpart is the following.
Let us consider two samples ${\bf X}=(X_1,\cdots, X_n)$ and ${\bf Y}=(Y_1,\cdots, Y_n)$. The so--called empirical Pearson correlation coefficient is given by
\begin{equation}\label{pearson}
r_n=\frac{\sum_{i=1}^n(X_i-\bar X_n)(Y_i-\bar Y_n)}{\sqrt{\sum_{i=1}^n(X_i-\bar X_n)^2 \sum_{i=1}^n(Y_i-\bar Y_n)^2}}\,,
\end{equation}
where $\bar X_n=\frac{1}{n}\sum_{k=1}^nX_k$ and $\bar Y_n=\frac{1}{n}\sum_{k=1}^nY_k$ are the empirical means of the samples.
Whenever $E({\bf X})$ and $E({\bf Y})$ are both known, we consider $\tilde r_n$:
\begin{equation}\label{pearsonknown}
\tilde r_n=\frac{\sum_{i=1}^n(X_i- E(X_i))(Y_i- E(Y_i))}{\sqrt{\sum_{i=1}^n(X_i- E(X_i))^2 \sum_{i=1}^n(Y_i- E(Y_i))^2}}\,.
\end{equation}

The study of the correlation coefficients is detailed in many references (see e.g.~\cite{Muir} or~\cite{Lot}) and it is shown that many competing€ correlation indexes are special cases of Pearson's correlation coefficient (\cite{RodNic}). 
The asymptotic behaviour of $(r_n)_n$, $(\tilde r_n)_n$ is worth considering. It is clear that when $\bf X$ and $\bf Y$ are independent, $r_n$, $\tilde r_n\to 0$ when $n\to \infty$. Moreover, whenever $({\bf X},{\bf Y})$ are sampled from a known distribution $(X,Y)$, $r_n$, $\tilde r_n\to \rho$ when $n\to \infty$. In this paper, we study Sharp Large Deviations (SLD) associated to these asymptotics.
\par
Large Deviations for empirical correlation coefficients have been studied by Si \cite{sishen} in the Gaussian case. We extend his results to SLD in the spherical and Gaussian cases. It can be noticed here that for the Gaussian case, we prove SLD on a resctricted domain of $\rho$ since the convexity properties of the functions are only true for $0\leq |\rho|\leq\rho_0$, where $\rho_0$ is explicitely defined. 
This point was not noticed in \cite{sishen} since the large deviations are given through a contraction principle which is actually not valid (the function used is not continuous). We stress the fact that things have to be handled in a different way for the case $|\rho|>\rho_0$, and there is no proof that the rate function should be the same.
\par
We consider here the asymptotic development of
\[P(r_n\geq c) \mbox{ or, equivalently } P(\tilde r_n\geq c)\]
for $0<c<1$. We follow the scheme of Bercu et al. \cite{BGL,BR} and split:
\[P(r_n\geq c)=A_nB_n\,,\]
where
\begin{eqnarray}
A_n=\exp[n(L_n(\lambda_c)-c\lambda_c)]\,,\\
B_n=\mathbb E_n[\exp[-n\lambda_c(r_n-c)]\1_{r_n\geq c}]\
\,,
\end{eqnarray}
$L_n$ is the normalized cumulant generating function (n.c.g.f.) of $r_n$, $L$ its limit as $n\to\infty$ and $\lambda_c$ is the unique $\lambda$ such that $L'(\lambda_c)=c$. We perform the following change of probability
\[\frac{dQ_n}{dP}=e^{\lambda_cnr_n-nL_n(\lambda_c)}\,,\]
and $\mathbb E_n$ is the expectation under this new probability $Q_n$. The key point here is to develop the characteristic function $\Phi_n$ of $\frac{\sqrt{n}(r_n-c)}{\sigma_c}$. We use an expansion already computed in the i.i.d. case by Cram\'er (see \cite{Cra}, Lemma 2, p.72) and Esseen \cite{Esseen}.
\par
Such studies have been done in the context of small ball deviations and Gaussian semi-norms by Ibragimov \cite{IbraLapla}, Li \cite{Li1}, Sytaya \cite{sytaya}, Zolotarev \cite{zolot}, with an anaytic point of vue and different asymptotics. Independently, with large deviations techniques, was done a similar work by Dembo, Meyer-Wolf and Zeitouni \cite{dembmay,mayerzeit}. We can also cite works in other contexts: Ben Arous \cite{BALapla} on asymptotic expansion of the heat kernel associated with an hypoelliptic operator (in small time), and Bolthausen \cite{bolthLapla} on the limiting behaviour of the partition function for random vectors in Banach spaces in a general i.i.d. case.
\par
The paper is organized as follows: in Sections \ref{sectspher} and \ref{sectgauss}, we present the SLD results in the spherical and Gaussian cases; Section \ref{proofs} is devoted to the proofs. In Section \ref{further}, we briefly extend our results to any order developments as we present an application to Bahadur exact slopes for the test based on $r_n$ in the Gaussian case. Finally, in an Appendix, we give some more details and references on the Laplace method.
\section{Spherical distribution}\label{sectspher}
In this section, we study empirical correlations coefficients \eqref{pearson} and \eqref{pearsonknown} under the following spherical distribution assumption. We denote by ${\bf v}'$ the transpose of vector $\bf v$.
\begin{assum}\label{assumspher}
Let ${\bf X}=(X_1,\cdots, X_n)'$ and ${\bf Y}=(Y_1,\cdots, Y_n)'$, $n\in\mathbb N$, $n>2$, be two independent random vectors where $\bf X$ has a $n$-variate spherical distribution with $P({\bf X}=(0,\cdots,0)')=0$ and ${\bf Y}$ has any distribution with $P({\bf Y}\in\{{\bf 1}\})=0$ (where ${\bf 1}=\{k(1,\cdots,1)'\,,k\in\mathbb R\}$).
\end{assum}  

\subsection{SLDP for $r_n$}
In order to derive SLD for $(r_n)_n$ we compute the n.c.g.f.
\begin{equation}\label{fonctionLn}
L_n(\lambda)=\frac{1}{n}\log E(e^{n\lambda r_n})\,.
\end{equation}
The asymptotics of $L_n$ are given in the following proposition:

\begin{prop}\label{loglapspher}
For any $\lambda\in\mathbb R$, we have
\begin{equation}\label{methlap}
E \left( e^{n\lambda r_n} \right)=\frac{\Gamma(\frac{n-1}{2})}{\pi^{1/2}\Gamma(\frac{n-2}{2})} e^{nh(r_0(\lambda))}\left(\frac{c_0(\lambda)}{\sqrt n}+\bigO{\frac{1}{n^{3/2}}}\right)\,,
\end{equation}
where
\begin{itemize}
\item
$h(r)=\lambda r+\frac{1}{2}\log(1-r^2)$,
\item 
$r_0(\lambda)$ is the unique root in $]-1,1[$ of $h'(r)=0$, i.e.
\begin{equation}\label{r0}
r_0(\lambda)=\frac{-1+\sqrt{1+4\lambda^2}}{2\lambda}\,,
\end{equation}
\item
$g(r)=(1-r^2)^{-2}$ and
$
c_0(\lambda)=\sqrt{\dfrac{2\pi}{|h''(r_0(\lambda))|}}g(r_0(\lambda))\,.
$
\end{itemize}
Therefore
\begin{equation}\label{loglap}
L_n(\lambda)=L(\lambda)-\frac{1}{n}\left[\frac{1}{2}\log\sqrt{1+4\lambda^2}-\frac{3}{2}\log\frac{1+\sqrt{1+4\lambda^2}}{2}\right]+\bigO{\frac{1}{n^2}}\,.
\end{equation}
where $L$ 
is the limit normalized $\log$--Laplace transform of $r_n$:
\begin{equation}
L(\lambda)=h(r_0(\lambda)).
\end{equation}
\end{prop}
The proof of this proposition is postponed to Section \ref{proofs}.
Now we have the following SLDP:
\begin{theo}\label{theospherical1}
For any $0<c<1$, under Assumption \eqref{assumspher}, we have
 \begin{equation}\label{mainspher}
 P(r_n\geq c)=\frac{e^{-nL^*(c)-\frac{1}{2}\log(1+4\lambda_c^2)+\frac{3}{2}\log\frac{1+\sqrt{1+4\lambda_c^2}}{2}}}{\lambda_c\sigma_c\sqrt{2\pi n}}(1+\smallO{1}),
  \end{equation}
where 
\begin{itemize}
\item $\lambda_c$ is the unique solution of 
$L'(\lambda)=c\,$, i.e. $\lambda_c=\dfrac{c}{1-c^2}$,
\item
$\sigma_c^2=L''(\lambda_c)			 =\dfrac{(1-c^2)^2}{1+c^2}$,
\item 
$L^*(y)=-\frac{1}{2}\log(1-y^2)\,.$
\end{itemize}
\end{theo}
\par
	{\it Proof:}
\par
To prove the SLD for $(r_n)_n$, we proceed as in Bercu et al.~\cite{BGL,BR}. 
The following lemma, which proof is given in the Section \ref{proofs},  gives some basic properties of $L$:

\begin{lem}\label{propriL}
Let $L(\lambda)=h(r_0(\lambda))$ where $h$ and $r_0$ are defined in Proposition \ref{loglapspher}, we have
\begin{itemize}
\item
$L$ is defined on $\mathbb R$ and $L$ is ${\cal C}^{\infty}$ on its domain.
\item $L$ is a strictly convex function on $\mathbb R$, $L$ reaches its minimum at $\lambda=0$. Moreover for any $\lambda\in\mathbb R$, $L'(\lambda) \in ]-1,1[$.
\item The Legendre dual of $L$ is defined on $]-1,1[$ and computed as
\begin{equation}
L^*(y)=\sup_{\lambda\in \mathbb R}\{\lambda y-L(\lambda)\}=-\frac{1}{2}\log(1-y^2) \, .
\end{equation}.
\end{itemize}
\end{lem}

Let $0<c<1$ and $\lambda_c>0$ such that $L'(\lambda_c)=c$.
Then
 \[L^*(c)=c\lambda_c-L(\lambda_c)\,,\]
 We denote by $\sigma_c^2=L''(\lambda_c)$, and define the following change of probability:
 \begin{equation}
 \frac{dQ_n}{dP}=e^{\lambda_cnr_n-nL_n(\lambda_c)} \, .
 \end{equation}
 The expectation under $Q_n$ is denoted by $E_n$.
We write
 \begin{equation}
 P(r_n\geq c)=A_nB_n\,,
 \end{equation}
 where
 \[A_n=\exp[n(L_n(\lambda_c)-c\lambda_c)]\,,\]
 \[B_n=E_n(\exp[-n\lambda_c(r_n-c)]\1_{r_n\geq c})\,.\]
 On the one hand, from \eqref{loglap}
 \[A_n=\exp[-nL^*(c)-\frac{1}{4}\log(1+4\lambda_c^2)+\frac{3}{2}\log\frac{1+\sqrt{1+4\lambda^2}}{2}] \left( 1+\bigO{\frac{1}{n}} \right) . \]
 On the other hand, let us denote by
 
 \[U_n=\frac{\sqrt n(r_n-c)}{\sigma_c},\]
 \[\Phi_n(u)=E_n(e^{iuU_n})=\exp(-\frac{iu\sqrt n}{\sigma_c}c+nL_n(\lambda_c+\frac{iu}{\sigma_c\sqrt n})-nL_n(\lambda_c))\,.\]
 We have the following technical results on $\Phi_n$, proved in Section \ref{proofs}.
 \begin{lem}\label{lemIPP}
 For any $K\in \mathbb N^*$, $\eta>0$, for $n$ large enough and any $u\in\mathbb R$,
 
 \begin{equation}\label{boundPhi}
 \displaystyle |\Phi_n(u)|\leq \frac{1}{|\lambda_c+\frac{iu}{\sigma_c\sqrt n}|^K}\frac{c_0^K(\lambda)}{c_0(\lambda)}(1+\eta)\,.
 \end{equation}
 where $c_0$ and $c_0^K$ are the first coefficients in Laplace's method (see Theorem \ref{theo_Lap}), corresponding respectively to 
 \[g(r)=(1-r^2)^{-2}\] and
 \[g^K(r)=(2r)^K(1-r^2)^{-K-2}\,.\]
 \end{lem}

From lemma above, choosing $K\geq 2$, we see that $\Phi_n$ is in $L^2$ and by Parseval formula,
 \begin{eqnarray*}
 B_n=E_n[e^{-\lambda_c\sigma_c\sqrt nU_n}\1_{U_n\geq 0}]
 =\frac{1}{2\pi}\int_{\mathbb R}\left(\frac{1}{\lambda_c\sigma_c\sqrt n+iu}\right)\Phi_n(u)du=\frac{C_n}{\lambda_c\sigma_c\sqrt{2\pi n}},
 \end{eqnarray*}
where
\[C_n=\frac{1}{\sqrt{2\pi}}\int_{\mathbb R}\left(1+\frac{iu}{\lambda_c\sigma_c\sqrt n}\right)^{-1}\Phi_n(u)du.\]
The key point here is to study the asymptotics of $\Phi_n$.

\begin{lem}\label{compPhi}
We have
\[\lim_{n\to\infty}\Phi_n(u)=e^{-u^2/2}\mbox{ and }\lim_{n\to\infty}C_n=1 \, .\]
\end{lem}

From lemma above, which proof is postponed to Section \ref{proofs}, we have equation \eqref{mainspher}. 

\hfill $\square$
 
 \subsection{Known expectation}

 In case $E(\bf X)$ and $E(\bf Y)$ are both known, we consider $\tilde r_n$ given in formula \eqref{pearsonknown} which can be written as follows
 \begin{equation}\label{pearsonknown_spher}
\tilde r_n=\frac{({\bf X}-E({\bf X}))' \, ({\bf Y}-E({\bf Y}))}{\Vert {\bf X}-E({\bf X})  \Vert \; \Vert {\bf Y}-E({\bf Y})  \Vert}\,.
 \end{equation}
We can derive a SLD result similar to the previous one. The following proposition gives the expression of the n.c.g.f. of $\tilde r_n$:
 \begin{prop}\label{loglapspher2}
For any $\lambda\in\mathbb R$, we have 
 \begin{equation}\label{methlap}
E(e^{n\lambda \tilde  r_n})=\frac{\Gamma(\frac{n}{2})}{\pi^{1/2}\Gamma(\frac{n-1}{2})} e^{nh(r_0(\lambda))}\left(\frac{\tilde c_0(\lambda)}{\sqrt n}+\bigO{\frac{1}{n^{3/2}}}\right)\,,
\end{equation}
where
 \begin{itemize}
\item
$h(r)=\lambda r+\frac{1}{2}\log(1-r^2)$,
\item 
$r_0(\lambda)$ is the unique root in $]-1,1[$ of $h'(r)=0$, i.e.
\[r_0(\lambda)=\frac{-1+\sqrt{1+4\lambda^2}}{2\lambda}\,,\]
\item
$\tilde  g(r)=(1-r^2)^{-3/2}$ and
$
\tilde c_0(\lambda)=\sqrt{\dfrac{2\pi}{|h''(r_0(\lambda))|}}\tilde  g(r_0(\lambda))\,.
$
\end{itemize}
The n.c.g.f. of $\tilde r_n$ is
 \begin{equation}\label{loglap2}
\tilde L_n(\lambda)=h(r_0(\lambda))-\frac{1}{n}\left[\frac{1}{2}\log\sqrt{1+4\lambda^2}-\log\frac{1+\sqrt{1+4\lambda^2}}{2}\right]+\bigO{\frac{1}{n^2}}\,.
\end{equation}
 \end{prop}
 This proposition is proved in Section \ref{proofs}. We have the following SLDP:

\begin{theo}\label{theospherical2}
For any $0<c<1$, under Assumption \eqref{assumspher}, we have
 \begin{equation}\label{mainspher2}
 P(\tilde r_n\geq c)=\frac{\exp^{-nL^*(c)-\frac{1}{4}\log(1+4\lambda_c^2)+\log\frac{1+\sqrt{1+4\lambda_c^2}}{2}}}{\lambda_c\sigma_c\sqrt{2\pi n}}(1+\smallO{1}).
  \end{equation}

\end{theo}
\par{\it Proof:}\par
 The proof of Theorem \ref{theospherical2} is exactly similar to the one of Theorem \ref{theospherical1} and formula \eqref{mainspher} is changed to \eqref{mainspher2} according to the way formula \eqref{loglap} is changed to \eqref{loglap2}.
 
 \hfill$\square$
 \section{Gaussian case}\label{sectgauss}

\begin{assum}\label{gauss}
Let $(X,Y)$ be a $\mathbb R^2$-valued Gaussian random vector where $\sigma_1^2=\mbox{Var}(X)$, $\sigma_2^2=\mbox{Var}(Y)$ and $\rho$ is the correlation coefficient: $\mbox{Cov}(X,Y)=\rho \sigma_1\sigma_2$.
We consider $({\bf X,\bf Y})=\{(X_i,Y_i), i=1,\cdots n\}$  an i.i.d. sample of $(X,Y)$.
\end{assum}
\subsection{General case}
We deal with the Pearson coefficient given in \eqref{pearson}.  As presviously mentioned,  Large deviations for $(r_n)_n$ are detailed in the paper of Si~\cite{sishen}. It can be noted that the contraction principle used by Si is not valid  here. The rate function is correct however, but only on some domain of $\rho$. We can give an expression of the normalized $\log$--Laplace transform $L_n$ given by \eqref{fonctionLn}.
 \begin{prop}\label{Gaussloglap}
 Let us define
 \[\rho_0:= \frac{\sqrt{3+2\sqrt{3}}}{3}\,.\]
For any $\lambda \in\mathbb R$ and $\rho$ such that $|\rho|\leq \rho_0$, we have the n.c.g.f. of $r_n$:
 \begin{equation}\label{loglap3}
 L_n(\lambda)= \overline{  h}(r_0(\lambda))  +\frac{1}{2}\log(1-\rho^2)+\frac{1}{n} \left[ \log \overline{  g}_{\rho}(r_0(\lambda))   -\frac{1}{2}\log |\overline{  h}''(r_0(\lambda))| \right] + \bigO{\frac{1}{n^2}} ,
 \end{equation}
  in which
   \begin{itemize}
\item
$\overline{  h}(r)=\lambda r-\log(1-\rho r)+\frac{1}{2}\log(1-r^2)$, 
\item 
$r_0(\lambda)$ is the unique real root in $]-1,1[$ of $\overline{  h}'(r)=0$,
\item
$\overline{  g}_{\rho}(r)= (1-\rho^2)^{-1/2}   (1-\rho r)^{3/2}(1-r^2)^{-2}$.
\end{itemize}
 \end{prop}
The proof of this proposition is postponed to Section \ref{proofs}.
We prove the following SLDP:
 \begin{theo}\label{Gauss_theo}
 For any $0 \leq\rho<c<1$ and  $|\rho| \leq \rho_0$ (with the notations of Proposition \ref{Gaussloglap}), we have
 \begin{equation}\label{mainGauss}
 P(r_n\geq c)=\frac{e^{-nL^*(c)+\log \overline{  g}_{\rho}(r_0(\lambda_c))   -\frac{1}{2}\log |\overline{  h}''(r_0(\lambda_c))|}}{\lambda_c\sigma_c\sqrt{2\pi n}}(1+\smallO{1}) \, ,
  \end{equation}
 where for any $-1<y<1$,
 \begin{equation}\label{gaussrate}
 L^*(y)=\log\left(\frac{1-\rho y}{\sqrt{(1-\rho^2)}\sqrt{(1-y^2)}}\right)\,.
 \end{equation}
 \end{theo}
 \par{\it Proof:}\par
Following the Proof of Theorem \ref{theospherical1}, we can easy obtain \eqref{mainGauss}. Note that the rate function in Si~\cite{sishen} matches our \eqref{gaussrate}. 

\hfill$\square$

\subsection{Known expectations}

In case $E(X)$ and $E(Y)$ are both known; and $\rho=0$, we have the following result

\begin{prop}\label{Gaussloglap2}
The n.c.g.f. of $\tilde r_n$ is given for any $\lambda\in\mathbb R$ by
\begin{equation}
L_n(\lambda)=h(u_0(\lambda))-\frac{1}{4n}\log(1+4\lambda^2)+\bigO{\frac{1}{n^2}},
\end{equation}
where 
\begin{itemize}
\item
$h(r)=\lambda r+\frac{1}{2}\log(1-r^2)$,
\item 
$u_0(\lambda)$ is the unique solution of $h'(\lambda)=0$ in $]-1,1[$.
\end{itemize}
\end{prop}

The proof is postponed to Section \ref{proofs}. The SLDP is therefore:

\begin{theo}
When $\rho=0$ and under Assumption \ref{gauss}, for $0<c<1$, we have
\begin{equation}
P(\tilde r_n\geq c)=\frac{e^{-nL^*(c)-\frac{1}{4}\log(1-4\lambda_c^2)}}{\lambda_c\sigma_c\sqrt n}(1+\smallO {1}),
\end{equation}
where $L^*$ is given in Theorem \ref{theospherical1}.
\end{theo}


 \section{Proofs}\label{proofs}
 \subsection{Proof of Proposition \ref{loglapspher}}
 We know from Muirhead (Theorem 5.1.1,~\cite{Muir}) that 
\[(n-2)^{1/2}\frac{r_n}{(1-(r_n)^2)^{1/2}}\]
has a Student's $t_{n-2}$-distribution. Hence the density function of $r_n$ is 

\begin{equation}\label{denspher}
f_n(r)=\frac{\Gamma(\frac{n-1}{2})}{\pi^{1/2}\Gamma(\frac{n-2}{2})}(1-r^2)^{(n-4)/2}\quad (-1<r<1).
\end{equation}
Applying Theorem \ref{theo_Lap}, we get
\begin{align*}
E \left( e^{n\lambda r_n}\right)	&=\int_{-1}^1e^{n\lambda r}f_n(r)dr=\int_{-1}^1e^{n\lambda r}\frac{\Gamma(\frac{n-1}{2})}{\pi^{1/2}\Gamma(\frac{n-2}{2})}(1-r^2)^{(n-4)/2}dr\\
		&=\frac{\Gamma(\frac{n-1}{2})}{\pi^{1/2}\Gamma(\frac{n-2}{2})}e^{nh(r_0(\lambda))}\left(\frac{c_0(\lambda)}{\sqrt n}+\bigO{\frac{1}{n^{3/2}}}\right).
\end{align*}
 where $h$, $r_0$ and $c_0$ are given in Proposition \ref{loglapspher}.

So we have
\begin{align}
E \left( e^{n\lambda r_n}\right)&=\frac{\Gamma(\frac{n-1}{2})}{\pi^{1/2}\Gamma(\frac{n-2}{2})}\sqrt{\frac{2\pi}{n}}e^{nh(r_0(\lambda))}\frac{g(r_0(\lambda))}{\sqrt{|h''(r_0(\lambda))|}}\left(1+\bigO{\frac{1}{n}}\right)\\\label{formula1}
&=\frac{\Gamma(\frac{n-1}{2})}{\Gamma(\frac{n-2}{2})}\sqrt{\frac{2}{n} } \, e^{nh(r_0(\lambda))}\frac{1}{(1-r_0(\lambda)^2)\sqrt{1+r_0(\lambda)^2}}\left(1+\bigO{\frac{1}{n}}\right)
\end{align}
From the duplication formula (see e.g. Olver \cite{Olv})
\[2^{2z-1}\Gamma(z)\Gamma(z+\frac{1}{2})=\sqrt \pi\Gamma(2z)\,,\]
as well as the Stirling formula (see \cite{Olv})
\[\log\Gamma(z)=z\log z-z-\frac{1}{2}\log z+\log\sqrt{2\pi}+\bigO{\frac{1}{Re(z)}}\,,\mbox{ as }Re(z)\to\infty\,,\]
formula \eqref{formula1} above becomes
\[E(e^{n\lambda r})=e^{nh(r_0(\lambda))}\frac{1}{(1-r_0(\lambda)^2)\sqrt{1+r_0(\lambda)^2}}\left(1+\bigO{\frac{1}{n}}\right)\,.\]
With the expression of $r_0$, we get formula \eqref{loglap}.
\subsection{Proof of Lemma \ref{propriL}}
We can explicit the full expression of $L$:
\begin{equation}
L(\lambda)=\frac{-1+\sqrt{1+4\lambda^2}}{2}-\frac{1}{2}\log(\frac{1+\sqrt{1+4\lambda^2}}{4})\,.
\end{equation}
It is easy to see that $L$ is defined on $\mathbb R$, ${\cal C}^{\infty}$ on its domain.
\par
From the definition of $L$ we can deduce
\begin{equation}
L'(\lambda)=r_0(\lambda)+h'(r_0(\lambda))=r_0(\lambda),
\end{equation}
and by construction of $r_0$, $L'\in]-1,1[$. Now we can compute
\begin{equation}
L''(\lambda)=r_0'(\lambda)=\frac{1}{2\lambda^2}\left(1-\frac{1}{\sqrt{1+4\lambda^2}}\right),
\end{equation}
and it is easily seen that $L''(\lambda)>0$ for any $\lambda\in\mathbb R^*$ and $L''(0)$ can be defined by continuity as $1$. Hence $L$ is strictly convex on $\mathbb R$ and has its minimum at $\lambda=0$. Moreover, if we have
\begin{equation*}
L'	(\lambda_c)		=	r_0(\lambda_c)=c,
\end{equation*}
then $0<c<1$ implies $\lambda_c>0$ and we can obtain 
\begin{equation*}
4 \lambda_c (\lambda_c(1-c^2)			-c)=0.
\end{equation*}
This leads us to the expression
\begin{equation*}
 \lambda_c=	\dfrac{c}{1-c^2}	\,.
\end{equation*}
Hence the preceding expression yields 
\begin{equation*}
 \sigma^2_c= L''( \lambda_c)=	\dfrac{(1-c^2)	^2}{1+c^2}	\,.
\end{equation*}
\subsection{Proof of Lemmas \ref{lemIPP} and \ref{compPhi}}

The proof of Lemma \ref{lemIPP} is based on iterated integrations by parts. We detail below the steps.
\begin{align*}
\Phi_n(u)		&=E_n(e^{iuU_n})=\int_{\mathbb R}e^{iu\frac{\sqrt{n}(r-c)}{\sigma_c}}f_n(r)e^{\lambda_cnr-nL_n(\lambda_c)}dr\\
&={\bf \Gamma}_n\,e^{-iu\frac{\sqrt{n}c}{\sigma_c}}e^{-nL_n(\lambda_c)}\int_{-1}^1e^{(iu\frac{\sqrt n}{\sigma_c}+\lambda_cn)r}(1-r^2)^{n/2-2}dr,
\end{align*}
where, for seek of simplicity, we denote by
\begin{equation}\label{Gammabf}
{\bf\Gamma}_n=\frac{\Gamma(\frac{n-1}{2})}{\pi^{1/2}\Gamma(\frac{n-2}{2})}\,.
\end{equation}
For $K\in\mathbb N^*$, performing $K$ integrations by part, since $f_n$ is zero at $-1$ and $1$ when $n$ is large enough, we get:
\begin{multline*}
\Phi_n(u)={\bf \Gamma}_n e^{-iu\frac{\sqrt{n}c}{\sigma_c}} e^{-nL_n(\lambda_c)}\times \cdots \\
\cdots\times\frac{(\frac{n}{2}-2)(\frac{n}{2}-3)\cdots(\frac{n}{2}-K-1)}{\left (iu\frac{\sqrt{n}}{\sigma_c}+\lambda_cn\right)^K}
\di\int_{-1}^1e^{(iu\frac{\sqrt n}{\sigma_c}+\lambda_cn)r}(-2r)^K(1-r^2)^{n/2-2-K}dr.
\end{multline*}
Hence,
\begin{equation*}
|\Phi_n(u)|\leq {\bf \Gamma}_n e^{-nL_n(\lambda_c)}
\frac{(\di\frac{n}{2}-2)(\frac{n}{2}-3)\cdots(\frac{n}{2}-K-1)}{\left |iu\frac{\sqrt{n}}{\sigma_c}+\lambda_cn\right|^K}
\int_{-1}^1e^{\lambda_cnr}(2r)^K(1-r^2)^{n/2-2-K}dr.
\end{equation*}
Using Laplace's method once again (see the Appendix), for a given $\eta>0$ we can find $N$ large enough such that for any $n\geq N$,
\begin{equation}
|\Phi_n(u)|\leq \frac{1}{|\lambda_c+\frac{iu}{\sqrt{n}\sigma_c}|^K}\frac{c_0^K(\lambda)}{c_0(\lambda)}(1+\eta)\,.
\end{equation}
\hfill$\square$
\par
To prove Lemma \ref{compPhi}, we first split $C_n$ into two terms:
\begin{equation}\label{split}
C_n=\frac{1}{\sqrt{2\pi}}\int_{|u|\leq n^{\alpha}}\left(1+\frac{iu}{\lambda_c\sigma_c\sqrt n}\right)^{-1}\Phi_n(u)du+\frac{1}{\sqrt{2\pi}}\int_{|u|> n^{\alpha}}\left(1+\frac{iu}{\lambda_c\sigma_c\sqrt n}\right)^{-1}\Phi_n(u)du.
\end{equation}
For the second term in the RHS of \eqref{split} we have
\begin{align*}
\left|\int_{|u|> n^{\alpha}}\frac{1}{\left(1+\frac{iu}{\lambda_c\sigma_c\sqrt n}\right)}\Phi_n(u)du\right|
&\leq \int_{|u|> n^{\alpha}}\frac{1}{\left|1+\frac{iu}{\lambda_c\sigma_c\sqrt n}\right|}\Phi_n(u)du\\
&\leq \int_{|u|> n^{\alpha}}\frac{1}{|\lambda_c|^K\left|1+\frac{iu}{\lambda_c\sigma_c\sqrt n}\right|^{K+1}}du\frac{c_0^K(\lambda_c)}{c_0(\lambda_c)}(1+\eta)\\
&\leq \frac{c_0^K(\lambda_c)}{|\lambda_c|^Kc_0(\lambda_c)}(1+\eta)
\int_{|u|> n^{\alpha}}\frac{1}{\left(1+\frac{u^2}{\lambda_c^2\sigma_c^2 n}\right)^{(K+1)/2}}du\\
&\leq \frac{c_0^K(\lambda_c)}{|\lambda_c|^Kc_0(\lambda_c)}(1+\eta)(\lambda_c^2\sigma_c^2 n)^{(K+1)/2}2\frac{n^{-\alpha K}}{K}.
\end{align*}

In order to have a negligible term, it is enough to have 
$-K\alpha+\frac{K+1}{2}<0$, i.e. fixing $K=3$, $\alpha=\frac{3}{4}$.
Now for the domain $\{|u|\leq n^{\alpha}\}$, we study more precisely the expression
\begin{equation}\label{Phi_n}
\Phi_n(u)=E_n(e^{iuU_n})=\exp\left[-\frac{iu\sqrt n}{\sigma_c}c+nL_n(\lambda_c+\frac{iu}{\sigma_c\sqrt n})-nL_n(\lambda_c)\right]\,.
\end{equation}
We first remark that $E(e^{n\lambda r_n})$ is analytic in $\lambda$ on $\mathbb R$, hence it can be expanded by analytic continuation and $L_n(\lambda+iy)$ for $\lambda,y\in\mathbb R$ is well defined.
From the analyticity, we can expand in Taylor series the expression \eqref{Phi_n} above.
\begin{align}\label{exprphi}
\Phi_n(\lambda_c)		&=\exp\{-iu\frac{\sqrt{n}c}{\sigma_c}+n\sum_{k=1}^{\infty}\left(\frac{iu}{\sigma_c\sqrt n}\right)^k\frac{L_n^{(k)}(\lambda_c)}{k!}\}	\nonumber  \\
&=\exp\{-iu\frac{\sqrt{n}c}{\sigma_c}+n\frac{iu}{\sigma_c\sqrt n}L_n'(\lambda_c)+n\sum_{k\geq 2}\left(\frac{iu}{\sigma_c\sqrt n}\right)^k\frac{L_n^{(k)}(\lambda_c)}{k!}\}.
\end{align}
We detail now a development of $L_n$ -- and its derivatives -- which will be useful in the whole paper.
 
\begin{techlem}\label{techlem1}
For any $\lambda\in\mathbb R$, we have
\begin{equation}\label{formtechlem}
L_n(\lambda)= h(r_0(\lambda))+\frac{1}{n}\log{\bf \Gamma}_n-\frac{1}{2n}\log n+\frac{1}{n}R_0(\lambda)+\frac{1}{n}\sum_{p\geq 1}\frac{R_p(\lambda)}{n^pp!},
\end{equation}
where ${\bf \Gamma}_n$ is defined in \eqref{Gammabf} and
\begin{eqnarray}\label{r0}
 R_0(\lambda)&=&\log c_0(\lambda),\\
 \label{rp}
R_p(\lambda)&=&\sum_{1\leq s\leq p}(-1)^{s-1}(s-1)!B_{p,s}(c_1,c_2,\cdots)c_0^{-s},
\end{eqnarray}
where the coefficients $c_i$ are given by  Laplace development (see Appendix) and $B_{p,s}$ are the partial exponential Bell polynomials (see \eqref{bellpart}).
\end{techlem}

{\it Proof of Technical Lemma \ref{techlem1}:}

From the Appendix we can develop
\begin{equation}
E(e^{n\lambda r_n})=\frac{\Gamma(\frac{n-1}{2})}{\pi^{1/2}\Gamma(\frac{n-2}{2})}\frac{e^{nh(r_0(\lambda)}}{\sqrt n}\sum_{p\geq 0}\frac{c_p(\lambda)}{(2p)!n^p},
\end{equation}
where
\begin{multline}\label{cp}
			c_p(\lambda)	=	\sqrt{\frac{2 \pi}{|	h''(r_0(\lambda))|}}	\,	\sum_{k=0}^{2p}	\left( \begin{array}{c} 2p \\ k \end{array} \right) 	g^{(2p-k)}(r_0(\lambda))		\\
	\cdot	\sum_{m=0}^k	B_{k,m} \left(	\frac{h^{(3)}(r_0(\lambda))	}{2.3}	, \dots		,\frac{h^{(k-m+3)}(r_0(\lambda))		}{(k-m+2)(k-m+3)}		\right)	\frac{(2m+2p-1)!!}{| h''(t_0)|^{m+p}}.
			\end{multline}

From Fa\`{a}  di Bruno formula (see e.g. formula [5c] of Comtet~\cite{comt1}):
\begin{eqnarray}
\log E(e^{n\lambda r_n})= nh(r_0(\lambda))+\log\left(\frac{\Gamma(\frac{n-1}{2})}{\sqrt n\pi^{1/2}\Gamma(\frac{n-2}{2})}\right)+\log {c_0(\lambda)}+\sum_{p\geq 1}\frac{R_p(\lambda)}{n^pp!},
\end{eqnarray}
where $R_p$ is defined in formula \eqref{rp} above. Hence the formula \eqref{formtechlem} is proven.
\hfill$\square$

\par
From expressions \eqref{rp} and \eqref{cp}, we see that $R_p$ is a polynomial in $g^{(s)}(r_0(\lambda))$ and $h^{(s)}(r_0(\lambda))$ where the derivatives are taken with respect to $r$.
The function $r_0(\lambda)$ is ${\cal C}^{\infty}$ on $\mathbb R$. We can therefore express the derivatives of $L_n$ as follows:
\begin{equation}\label{der}
L_n^{(k)}(\lambda)=L^{(k)}(\lambda)+\frac{R_0^{(k)}(\lambda)}{n}+\frac{1}{n}\sum_{p\geq 1}\frac{R_p^{(k)}(\lambda)}{n^pp!}.
\end{equation}

Back to formula \eqref{exprphi}, 
and from the choice of $\lambda_c$, we have
\[\frac{\partial}{\partial\lambda}h(r_0(\lambda))\Big|_{\lambda=\lambda_c}=L'(\lambda_c)=c\]
and
\begin{multline}
\Phi_n(u)=\exp\{\frac{iu\sqrt n}{\sigma_c}[L_n'(\lambda_c)-c]+n\sum_{k\geq 2}\left(\frac{iu}{\sigma_c\sqrt n}\right)^k\frac{L_n^{(k)}(\lambda_c)}{k!}\}\\
\label{devphi} =\exp\{\frac{iu}{\sqrt n\sigma_c}[R_0'(\lambda)+\sum_{p\geq 1}\frac{R_p'(\lambda)}{n^pp!}]-\frac{u^2}{2\sigma_c^2}L_n''(\lambda_c)+n\sum_{k\geq 3}\left(\frac{iu}{\sigma_c\sqrt n}\right)^k\frac{L_n^{(k)}(\lambda_c)}{k!}\}\\
=\exp\{
-\frac{u^2}{2}
+\sum_{k= 3}^{2p}\left( \frac{iu}{\sigma_c\sqrt n}\right)^k\frac{nL^{(k)}(\lambda_c)}{k!}
+\sum_{k= 1}^{2p}\left(\frac{iu}{\sigma_c\sqrt n}\right)^k\frac{R_0^{(k)}(\lambda_c)}{k!}+\sum_{k\geq 1}\left(\frac{iu}{\sigma_c\sqrt n}\right)^k\frac{1}{k!}\sum_{p\geq 1}\frac{R_p^{(k)}(\lambda_c)}{n^pp!}\}.
\end{multline}
For $p$ large enough such that $\{u^k/(\sqrt{n})^{k+2p}\}$ is bounded on $\{|u|\leq n^{\alpha}\}$, we can have a uniform bound on the rest of the sum in the last term on the RHS above.
Hence we can write, for a given $m\in\mathbb N$ large enough
\begin{multline}
\Phi_n(u)=\exp\{-\frac{u^2}{2}
+\sum_{k= 3}^{2m+3}\left( \frac{iu}{\sigma_c\sqrt n}\right)^k\frac{nL^{(k)}(\lambda_c)}{k!}
+\sum_{k= 1}^{2m+1}\left(\frac{iu}{\sigma_c\sqrt n}\right)^k\frac{R_0^{(k)}(\lambda_c)}{k!}\\
+\sum_{k= 1}^{2m+1}\sum_{p= 1}^{s(m)}\left(\frac{iu}{\sigma_c\sqrt n}\right)^k\frac{1}{k!}\frac{R_p^{(k)}(\lambda_c)}{n^pp!}\}+
O(\frac{1+|u|^{2m+4}}{n^{m+1}})
\}.
\end{multline}
We follow the scheme of Cramer~\cite{Cra} Lemma 2, p.72 (see also Bercu and Rouault~\cite{BR}), and we get the wanted results.

\hfill $\square$
\begin{rem}
A thorough study of expressions $L_n^{(k)}$ and $R_p^{(k)}$ are given in \cite{TienThesis}.
\end{rem}

 \subsection{Proof of Proposition \ref{loglapspher2}}
 By symmetry, the mean $E({\bf X})=0$ if it exists. Then, $\tilde r_n$ from \eqref{pearsonknown_spher} becomes
 \begin{equation}\label{pearsonknown2}
\tilde r_n= \frac{{\bf X}' ({\bf Y}-E({\bf Y}))}{\Vert {\bf X}  \Vert \; \Vert {\bf Y}-E({\bf Y})  \Vert}.
\end{equation}
 Applying Theorem 1.5.7 from Muirhead~\cite{Muir}, with $\alpha =\dfrac{{\bf Y}-E({\bf Y})}{\Vert {\bf Y}-E({\bf Y})  \Vert} \in \mathbb{R}^n$, then
\[(n-1)^{1/2}\frac{\tilde r_n}{(1-\tilde r_n^2)^{1/2}}\]
has a $t_{n-1}$-distribution. Comparing to $r_n$, the degree of the $t$-distribution is one degree less than $\tilde r_n$.

Hence the density function of $\tilde r_n$ is 
\begin{equation}\label{denspher2}
\frac{\Gamma(\frac{n}{2})}{\pi^{1/2}\Gamma(\frac{n-1}{2})}(1-r^2)^{(n-3)/2}\,,\quad (-1<r<1).
\end{equation}
Applying  Laplace's method we get
\begin{align*}
E \left(  e^{n\lambda \tilde r_n} \right)		&=\int_{-1}^1e^{n\lambda r}\frac{\Gamma(\frac{n}{2})}{\pi^{1/2}\Gamma(\frac{n-1}{2})}(1-r^2)^{(n-3)/2}dr\\
&=\frac{\Gamma(\frac{n}{2})}{\pi^{1/2}\Gamma(\frac{n-1}{2})} e^{nh(r_0(\lambda))}\left(\frac{\tilde c_0(\lambda)}{\sqrt n}+\bigO{\frac{1}{n^{3/2}}}\right),
\end{align*}
 where $h$, $r_0$ and $c_0$ are given in Proposition \ref{loglapspher2}. Then
 \begin{align}
 E \left(  e^{n\lambda \tilde r_n} \right)	&=\frac{\Gamma(\frac{n}{2})}{\pi^{1/2}\Gamma(\frac{n-1}{2})}\sqrt{\frac{2\pi}{n}}e^{nh(r_0(\lambda))}\frac{\tilde g(r_0(\lambda))}{\sqrt{|h''(r_0(\lambda))|}}\left(1+\bigO{\frac{1}{n}}\right) \nonumber \\\label{formula2}
&=e^{nh(r_0(\lambda))}\frac{1}{\sqrt{(1-r_0^2(\lambda))(1+r_0^2(\lambda))}}\left(1+\bigO{\frac{1}{n}}\right).
 \end{align}
And we can obtain formula \eqref{loglap2} from the expression of $r_0$.
\subsection{Proof of Proposition \ref{Gaussloglap}} \label{discuss Si}

From Muirhead, we know that the density function of a $n+1$ sample correlation coefficient $r_{n+1}$ is given by
\begin{eqnarray*}
\frac{(n-1)\Gamma(n)}{\Gamma(n+1/2)\sqrt{2\pi}}(1-\rho^2)^{n/2} (1-\rho r)^{-n+1/2}(1-r^2)^{(n-3)/2} \qquad\qquad\qquad
\\
{_2 F_1} \left(\frac{1}{2},\frac{1}{2};n+\frac{1}{2};\frac{1}{2}(1+\rho r)\right) \qquad (-1<r<1).
\end{eqnarray*} 
where $_2F_1$ is the hypergeometric function (see~\cite{Olv}).
Hence Laplace transform is
\begin{multline*}
E \left(e^{(n+1)\lambda r_{n+1}}\right)=\frac{(n-1)\Gamma(n)}{\Gamma(n+1/2)
\sqrt{2\pi}}(1-\rho^2)^{n/2}\\
\int_{-1}^1e^{(n+1)\lambda r}(1-\rho r)^{-n+1/2}(1-r^2)^{(n-3)/2}
{_2 F_1} \left(\frac{1}{2},\frac{1}{2};n+\frac{1}{2};\frac{1}{2}(1+\rho r)\right)dr.
\end{multline*}
Looking for a limit as $n\to \infty$, we can use the following result due to Temme~\cite{Temme1,Temme2} (see also~\cite{Hypergeo}): the function ${_2F_1}$ has the following Laplace transform representation
\begin{equation}
{_2F_1}(a,b,c;z)=\frac{\Gamma(c)}{\Gamma(b)\Gamma(c-b)}\int_0^1\frac{t^{b-1}(1-t)^{c-b-1}}{(1-zt)^a}dt
\end{equation}
and 
\begin{equation}
_2F_1(a,b,c+\lambda;z)\sim \frac{\Gamma(c+\lambda)}{\Gamma(c+\lambda-b)}\sum_{s=0}^{\infty}\textrm{f}_s(z)\frac{(b)_s}{\lambda^{b+s}}\, ,
\end{equation}
where the equivalent is for $\lambda\to+\infty$ and 
\[\textrm{f}(t)=\left(\frac{e^t-1}{t}\right)^{b-1}e^{(1-c)t}(1-z+ze^{-t})^{-a}\,,\]
\[\textrm{f}(t)=\sum_{s=0}^{\infty}\textrm{f}_s(t)t^s.\]
In our case, we get as $n\to\infty$:
\begin{equation}
_2F_1\left( \frac{1}{2},\frac{1}{2},\frac{1}{2}+n;\frac{1}{2}(1+\rho r) \right)\sim 
\frac{\Gamma(\frac{1}{2}+n)}{\Gamma(n)}		\left(\frac{1}{\sqrt n}+\frac{2+\rho r}{8n^{3/2}}+\smallO{\frac{1}{n^{3/2}}}		\right) \,.
\end{equation}
Hence we have to deal with the following integral:
\begin{equation}
\int_{-1}^1e^{(n+1)\lambda r}(1-\rho r)^{-n+1/2}(1-r^2)^{(n-3)/2}		\left( 1+\frac{2+\rho r}{8n}+\smallO{\frac{1}{n}}    \right)  dr.
\end{equation}
Neglecting the terms of lower order in $n$ we focus on
\begin{equation}
\int_{-1}^1e^{(n+1)\lambda r}(1-\rho r)^{-n+1/2}(1-r^2)^{(n-3)/2}dr=\int_{-1}^1e^{n \overline{  h} (r)} \overline{  g}(r)dr\,,
\end{equation}
where
\begin{equation}\label{barh}
\overline{  h}(r)=\lambda r-\log(1-\rho r)+\frac{1}{2}\log(1-r^2)\,,
\end{equation}
\[\overline{  g}(r)=e^{\lambda r}\sqrt{(1-\rho r)}(1-r^2)^{-3/2}\,.\]
The following lemma details the properties of the function $\bar h$:
\begin{lem}
For any $\rho\in ]-1,1[$ and $r\in ]-1,1[$, the function $\bar h$ of formula \eqref{barh} is defined for any $\lambda\in\mathbb R$.
Moreover the equation $\bar h'(r)=0$ has at least one solution in $]-1,1[$ and $\bar h''(r)<0$ on $]-1,1[$ for any $|\rho|\leq \rho_0$
where $\rho_0=\dfrac{\sqrt{3+2\sqrt{3}}}{3}\,.$

\end{lem}
\par{\it Proof:}\par
We compute easily
\[\overline h'(r)= \lambda+\dfrac{\rho}{1-\rho r}- \dfrac{r}{1-r^2}\]
and see that $H(r)=\bar h'(r)(1-r^2)=0$ has at least one root in $]-1,1[$ (since $H(-1)H(1)<0$). Hence there exists at least one solution $r_0 \in ]-1,1[$ such that $\bar h'(r)=0$.
Next, we compute
\[\bar h''(r)= \dfrac{\rho^2}{(1-\rho r)^2}-\dfrac{1+r^2}{(1-r^2)^2} \]
and we have 
\[\bar h''(r)<0 \mbox { for any }r\in ]-1,1[ \iff |\rho|\leq \rho_0:=\dfrac{\sqrt{3+2\sqrt{3}}}{3}\,.	 \qquad\qquad\square	\]

We know from Si~\cite{ sishen} that the rate function in this case is
\begin{equation}
I_{\rho}(s)=\log\left(\frac{1-\rho s}{\sqrt{(1-\rho^2)}\sqrt{(1-s^2)}}\right)\mbox{ for  }-1<s<1\,.
\end{equation}

As previously said, even if this function was obtained by a contraction principle which is not applicable here (the function involved is not continuous, see Dembo and Zeitouni for more details~\cite{DZ}), we claim that the expression of the rate function above is nevertheless correct in the given domain $\{|\rho|\leq \rho_0\}$. We prove it below.
We have
\[L(\lambda)=\bar h(r_0(\lambda))+\frac{1}{2}\log(1-\rho^2)\,,\]
where $r_0$ satisfies
\[\bar h'(r_0(\lambda))=0.\]
Now we compute
\begin{equation}
L'(\lambda)=r_0(\lambda)+r'_0(\lambda)\bar h'(r_0(\lambda))=r_0(\lambda).
\end{equation}
For every $-1<c<1$ and $\lambda_c$ such that $L'(\lambda_c)=c$, we have
\begin{align*}
L^*(c)	&=c\lambda_c-L(\lambda_c)\\
&=c\lambda_c-\{\lambda_cr_0(\lambda_c)+\frac{1}{2}\log(1-r_0^2(\lambda_c))-\log (1-\rho r_0(\lambda_c))+\frac{1}{2}\log(1-\rho^2)\}\\
&=-\frac{1}{2}\log(1-c^2)+\log (1-\rho c)-\frac{1}{2}\log(1-\rho^2)=\log\frac{1-\rho c}{\sqrt{1-c^2}\sqrt{1-\rho^2}}.
\end{align*}

From the dual properties of Legendre transform, the condition of Laplace's method $\overline{  h}''(r)<0$ is compatible with the condition of convexity of $I_{\rho}$ in $]-1,1[$.
Indeed, for $\rho_0<|\rho|<1$, $I_{\rho}$ is not convex. 
From that point, under condition $|\rho| \leq \rho_0$, we can get
\begin{align}
E \left(e^{(n+1)\lambda r_{n+1}}\right)&=\frac{n-1}{\sqrt{2n\pi}}(1-\rho^2)^{n/2} \sqrt{\frac{2\pi}{n}}e^{n\overline{  h}(r_0(\lambda))}\frac{\overline{  g}(r_0(\lambda))}{\sqrt{|\overline{  h}''(r_0(\lambda))|}}\left(1+\bigO{\frac{1}{n}}\right)\\\label{formula3}
&=e^{(n+1)\overline{  h}(r_0(\lambda))}\frac{(1-\rho^2)^{n/2} (1-\rho r_0(\lambda))^{3/2}}{(1-r_0^2(\lambda))^2 \sqrt{|\overline{  h}''(r_0(\lambda))|}} \left(1+\bigO{\frac{1}{n}}\right)\,.
\end{align}
We can adjust the size of sample into $n$ and obtain
\begin{equation}
E \left(e^{n\lambda r_{n}}\right)=e^{n\overline{  h}(r_0(\lambda))}\frac{(1-\rho^2)^{(n-1)/2} (1-\rho r_0(\lambda))^{3/2}}{(1-r_0^2(\lambda))^2		\sqrt{|\overline{  h}''(r_0(\lambda))|} }\left(1+\bigO{\frac{1}{n}}\right),
\end{equation}
which leads us to \eqref{loglap3}. We give below two graphics, one for $\rho=\rho_0-0.1$ and one for $\rho=\rho_0+0.1$. We can clearly see the change of convexity.
\medskip\pn
\begin{figure}[ht!]
\begin{center}
\includegraphics[width=.5\textwidth]{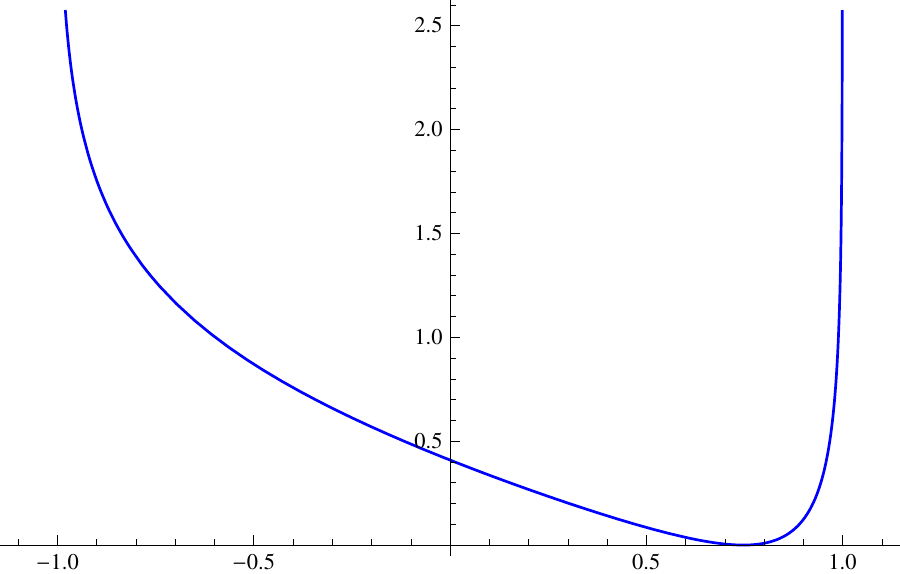}

\end{center}
\caption{$I_{\rho}$ for $\rho=\rho_0-0.1$}
\end{figure}
\begin{figure}[ht!]
\begin{center}
\includegraphics[width=.5\textwidth]{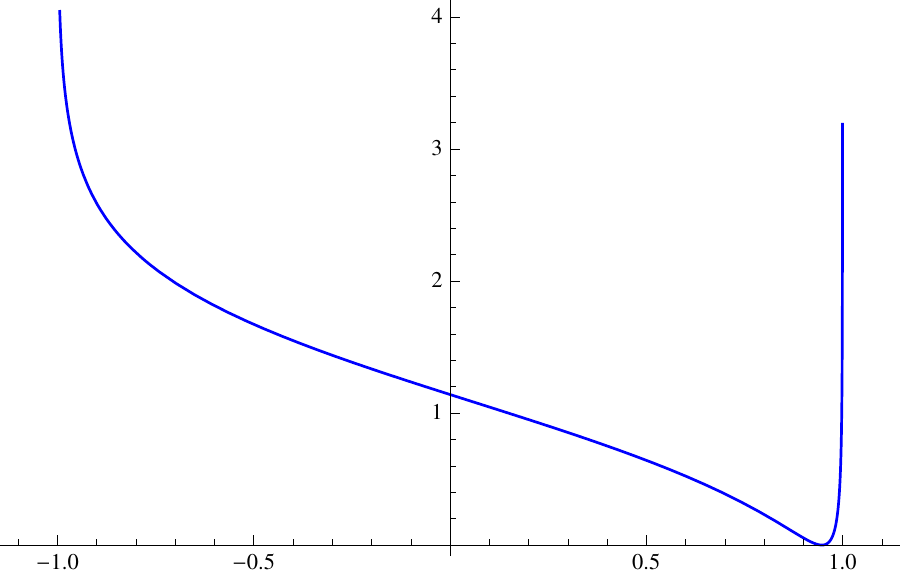}
\end{center}
\caption{$I_{\rho}$ for $\rho=\rho_0+0.1$}
\end{figure}

\subsection{Proof of Proposition \ref{Gaussloglap2}}
For the asymptotics of $L_n$ in this case, we follow the steps of Si~\cite{sishen}. Up to considering $X_1=X-E(X)$ and $Y_1=Y-E(Y)$, we can boil down to $E(X)=E(Y)=0$.

If we denote by $\langle,\rangle$ the Euclidean scalar product in $\mathbb R^2$, and
\[\tilde X=\left ( \frac{X_1}{\sqrt{\sum_{i=1}^nX_i^2}},\cdots,\frac{X_n}{\sqrt{\sum_{i=1}^nX_i^2}}\right)\,,\quad\tilde Y=\left (\frac{Y_1}{\sqrt{\sum_{i=1}^nY_i^2}},\cdots,\frac{Y_n}{\sqrt{\sum_{i=1}^nY_i^2}}\right)\,,\]
therefore
\begin{equation}
\tilde r_n=\langle \tilde X,\tilde Y\rangle.
\end{equation}
Large deviations for $(\tilde r_n)_n$ are proved in~\cite{sishen}. We derive here the corresponding sharp principle.
Since $\tilde X$, $\tilde Y$ are independent random variables with uniform distribution $\tilde \sigma_n$ on the unit sphere ${\cal S}^{n-1}$ of $\mathbb R^n$, we can compute
\begin{align}
E  \left( e^{\lambda \tilde r_n} \right)&=\iint\limits_{{\cal S}^{n-1}\times {\cal S}^{n-1}}e^{\lambda\langle x,y\rangle}\tilde\sigma_n(dx)\tilde\sigma_n(dy)dxdy\\
&=\frac{a_{n-1}}{a_n}\int_{-1}^{1}e^{\lambda u}\left(\sqrt{1-u^2}\right)^{n-1}du\,,
\end{align}
where $a_n$ is the area of the unit sphere:
\[a_i=\frac{2\pi ^{\frac{i+1}{2} }}{\Gamma(\frac{i+1}{2})}\,.\]
In order to get the SLD, we want to compute the normalized $\log$--Laplace transform: for any $\lambda\in\mathbb R$,
from Stirling formula (see \cite{Olv}), we get easily
\[\frac{a_{n-1}}{a_n}=\sqrt{\frac{n}{2\pi }}\left(1+\bigO{\frac{1}{n}}\right).\]
Then we can write
\[
\int_{-1}^1e^{n\lambda u}\left(\sqrt{1-u^2}\right)^{n-1}du=\int_{-1}^1e^{nh(u)}g(u)du\,,
\]
where $h(u)=\lambda u+\frac{1}{2}\log(1-u^2)$ and $g(u)=\frac{1}{\sqrt{1-u^2}}$. We apply Laplace's method to get:
\begin{equation}
\int_{-1}^1e^{nh(u)}du=e^{nh(u_0(\lambda))}\left(\frac{c_0(\lambda)}{\sqrt n}+\bigO{\frac{1}{n^{3/2}}}\right)\,,
\end{equation}
where 
\[u_0(\lambda)=\frac{-1+\sqrt{1+4\lambda^2}}{2\lambda}\,,\quad
 c_0(\lambda)=\sqrt{\frac{2\pi}{h''(u_0(\lambda))|}}g(u_0(\lambda))\,.\]
 This leads to
 \begin{equation}
 L_n(\lambda)=h(u_0(\lambda))-\frac{1}{2n}\log(1+4\lambda^2)+\bigO{\frac{1}{n^2}}.
 \end{equation}

\section{Further results}\label{further}
\subsection{Any order development}
We present in this section a way to extend the results of Sections \ref{sectspher} and \ref{sectgauss} to higher orders. Moreover, whenever functions involved are smooth enough, these techniques can be applied and the asymptotics are given in other cases.


\begin{theo}\label{g}
In the framework of Sections \ref{sectspher} and \ref{sectgauss}, for any $0<c<1$, there exists a sequence $(\delta_{c,k})_k$ such that
\begin{equation}
 P(r_n\geq c)=\frac{e^{-nL^*(c)+R_0(\lambda_c)}}{\lambda_c\sigma_c\sqrt{2\pi n}}		\left[ 1+ \sum_{k=1}^{p}  \frac{\delta_{c,k}}{n^k} +\bigO {\frac{1}{n^{p+1}}} \right].
  \end{equation}
\end{theo}
\par{\it Proof:}\par
For seek of simplicity, we only present here the proof for $(r_n)_n$ in the spherical case. Similarly to the proof of Theorem \ref{theospherical1}, we briefly give the main ideas:
From the decomposition $P(r_n\geq c)=A_nB_n$, in which
\begin{align*}
A_n		&=\exp[n(L_n(\lambda_c)-c\lambda_c)] \\
		&= 	\exp[-nL^*(c)+R_0(\lambda_c)+\sum_{p\geq 1}\frac{R_p(\lambda_c)}{n^p(2p)!} ]\\
		&= \exp[-nL^*(c)+R_0(\lambda_c))]\left( 1+\sum_{p\geq 1}\frac{\eta_p(\lambda_c)}{n^p(2p)!}\right)\,,
\end{align*}
where $(\eta_p)_p$ is a sequence of smooth functions of $\lambda$.
Recall that we can develop $L_n$ as in \eqref{formtechlem} and $L_n^{(k)}$ as in \eqref{der}
and from the development of $\Phi_n$ in \eqref{devphi},
\begin{equation}
\left(1+\frac{iu}{\lambda_c\sigma_c\sqrt n}\right)^{-1}\Phi_n(u)= e^{-\frac{u^2}{2 \sigma_c}}		\left( 1+ \sum_{k=1}^{2p+1} \frac{P_{p,k}(u)}{n^{k/2}}	+	\frac{1+u^{6(p+1)}}{n^{p+1}}\bigO {1}\right),
\end{equation}
where $P_{p,k}$ are polynomials in odd powers of $u$ for $k$ odd, and polynomials in even powers of $u$ for $k$ even.
From that points, we can complete the proof of Theorem \ref{g}.

\hfill$\square$
\subsection{Correlation test and Bahadur exact slope}

\subsubsection{Bahadur slope}
Let us recall here some basic facts about Bahadur exact slopes of test
statistics.  For a reference, see~\cite{Bah} and~\cite{Nik}.
Consider a sample $X_1, \cdots, X_n$ having common law $\mu_{\theta}$
depending on a parameter $\theta \in \Theta$.  To test
$(H_0): \theta \in \Theta_0$ against the alternative
$(H_1): \theta \in \Theta_1= \Theta \backslash
\Theta_0$, we use a test statistic $S_n$, large values of $S_n$
rejecting the null hypothesis.  The $p$-value of this test is by
definition $G_n(S_n)$, where
\begin{equation*}
  G_n(t) = \sup_{\theta\in\Theta_0} P_\theta(S_n \geq t).
\end{equation*}
The Bahadur exact slope $c(\theta)$ of $S_n$ is then given by the
following relation 
\begin{equation}
  \label{slope}
   c(\theta)=-2\liminf_{n\to\infty} \frac{1}{n} \log\left(G_n(S_n)\right).
\end{equation}
Quantitatively, for $\theta \in \Theta_1$, the larger $c(\theta)$ is,
the faster $S_n$ rejects $H_0$.

A theorem of Bahadur (Theorem 7.2  in~\cite {Baha}) gives the
following characterization of $c(\theta)$: suppose that $\lim_nn^{-1/2}S_n = b(\theta)$ for any $\theta\in\Theta_1$, and that 
$\lim_nn^{-1}\log\left(G_n(n^{1/2}t)\right) = -I(t)$ under any $\theta\in\Theta_0$.
If $I$ is continuous on an interval containing $b(\Theta_1)$,
then $c(\theta)$ is given by:
\begin{equation}
  \label{slope1}
  c(\theta) = 2 I(b(\theta)) \, .
\end{equation}

\subsubsection{Correlation in the Gaussian case}
In the Gaussian case, under Assumption \ref{gauss}, we have the following strong law of large numbers:
\begin{equation}\label{LLN}
r_n\to \rho=\mbox{cov}(X,Y)\,.
\end{equation}
We wish to test $H_0: \,\,\rho=0$ against the alternative $H_1:\,\, \rho\neq 0$.
It is obvious that under $H_1$,
\[\lim_{n\to\infty}r_n=\rho \,,\]
and this limit is continuous when $\rho\neq 0$.

Besides, we have here
\[G_n(t)=\sup_{\rho\in\Theta_0}P_{\rho}(\sqrt n r_n\geq t)\]
and 
\[\frac{1}{n}\log G_n(\sqrt nt)\to -\frac{1}{2}\log(1-t^2)\,.\]
Therefore the Bahadur slope is
\begin{equation}
c(\rho)=\log (1-\rho^2)\,.
\end{equation}
We show that this statistic is optimal in a certain sense.
In the framework above,  to test  $\theta\in \Theta_0$ against the alternative $\theta\in\Theta_1$ we define the likelihood ratio:
\[\lambda_n=\frac{\sup_{\theta\in\Theta_0}\prod_{i=1}^n\mu_{\theta}(x_i)}{\sup_{\theta\in\Theta_1}\prod_{i=1}^n\mu_{\theta}(x_i)}\]
and the related statistic:
\begin{equation}\label{loglik}
\hat S_n=\frac{1}{n}\log\lambda_n.
\end{equation}
Bahadur showed in~\cite{Bah2} that $\hat S_n$ is optimal in the following sense: for any $\theta\in\Theta_1$,
\begin{equation}
\lim_{n\to\infty}\frac{1}{n}\log G_n(\hat S_n)=-J(\theta),
\end{equation}
where $J$ is the infimum of the Kullback--Leibler information:
\begin{equation}
J(\theta)=\inf\{K(\theta,\theta_0),\theta_0\in\Theta_0\}
\end{equation}
and
\begin{equation}
K(\theta,\theta_0)=-\int\log[\frac{\mu_{\theta_0}(x)}{\mu_{\theta}(x)}]d\mu_{\theta}\,.
\end{equation}
\begin{defi}
Let $T_n$ be a statistic in the parametric framework defined above, then if $c(\theta)$ is the Bahadur slope of $T_n$, we have
\[c(\theta)\leq 2J(\theta)\]
and $T_n$ is said to be optimal if the upper bound is reached.
\end{defi}
We have the following result on the statistic $r_n$.
\begin{prop}
The sequence of empirical coefficients $(r_n)_n$ is asymptotically optimal in the Bahadur sense (\cite{Bah2}).
\end{prop} 
\par{\it Proof:}\par
We can easily compute the Kullback--Liebler information in this case:

 Let $\theta=(\mu,\Sigma)$ corresponds to the distribution of $(X,Y)$ in the case $\theta\in\Theta_1$ and $\theta=(\mu_0,\Sigma_0)$ for $\theta\in\Theta_0$. Since $\rho=0$ in the case $\theta\in\Theta_0$, the matrix $\Sigma_0$ is diagonal.
\begin{equation}\label{kulb}
K(\theta,\theta_0)=-\frac{1}{2}\log|\Sigma|+\frac{1}{2}\log|\Sigma_0|-1+\frac{1}{2}\mbox{tr}\Sigma_0^{-1}[\Sigma-(\mu-\mu_0)^t(\mu-\mu_0)],
\end{equation}
where $|\Sigma|$ stands for the determinant of $\Sigma$.
The infimum in \eqref{kulb} is reached when $\mu_0=\mu$ and the diagonal terms in $\Sigma_0$ are the ones of $\Sigma$.

Hence,
\[J(\theta)=\inf_{\theta_0\in\Theta_0}K(\theta,\theta_0)=-\frac{1}{2}\log|\Sigma|+\frac{1}{2}\log\sigma_{11}+\frac{1}{2}\log\sigma_{22}
=-\frac{1}{2}\log(1-\rho^2).\]
\hfill$\square$

\section*{Appendix: Laplace method}\label{appendix}
We present here some well known results about asymptotics of Laplace transforms. More precisely, we consider integrals of type
\begin{equation}\label{Lap}
I(x)=\int_a^be^{xp(t)}q(t)dt
\end{equation}
and its asymptotics as $x\to\infty$. Details and references can be found in Olver \cite{Olv} and Queffelec and Zuily \cite{queffzuil}. The explicit computations are also done in \cite{TienThesis}.
Let us first recall some definitions (for more details, see Comtet \cite{comt1,comt2}).
\begin{defi}\label{bellpartdef}
Partial exponential Bell polynomials are defined for any positive integers $k\leq n$ by
\begin{equation}\label{bellpart}
B_{n,k}(x_1,x_2,\cdots,x_{n-k+1})=\sum\frac{n!}{c_1!c_2!\cdots c_{n-k+1}!} \left(		\dfrac{x_1}{1!}\right)  ^{c_1}	\left(		\dfrac{x_2}{2!}\right)  ^{c_2}		\cdots \, \left(		\dfrac{x_{n-k+1}}{(n-k+1)!}\right)  ^{c_{n-k+1}}	,
\end{equation}
where the sum is taken over all positive integers $c_1,c_2\cdots,c_{n-k+1}$ such that
\begin{eqnarray*}
c_1+c_2+\cdots +c_{n-k+1}=k,\\
c_1+2c_2+\cdots +(n-k+1)c_{n-k+1}=n.
\end{eqnarray*}
\end{defi}
\begin{defi}\label{bellcompl}
The complete exponential Bell polynomials are defined by
\begin{eqnarray*}
B_0=1,\\
\forall n\geq 1\,,\quad B_n=\sum_{k=1}^nB_{n,k}.
\end{eqnarray*}
where $B_{n,k}$ are partial exponential Bell polynomials defined above.
\end{defi}

 \begin{theo} \label{theo_Lap}
                Let $(a,b)$ be a non-empty open interval, possibly non bounded and  $t_0$ be some point in $(a,b)$. Denote by $V_{t_0}$ a neighborhood of $t_0$ such that $p,q:(a,b)\to\mathbb R$ are functions of class $\mathcal{C}^{\infty}(V_{t_0})$.

		We suppose that
			\begin{itemize}
				\item[i)]  $p$ is measurable on $(a,b)$,
				\item[ii)] The maximum of $p$ is reached at $t_0$ (i.e. $p'(t_0)=0$ and $p''(t_0)<0$),
				\item[iii)]There exists $x_0$ such that $ \displaystyle \int_a^b	e^{x_0 p(t)} \, |q(t)| dt < + \infty$.
			\end{itemize}
		Then there exist coefficients $c_0(t_0), c_1(t_0), \dots $ depending on derivatives of $p$ and $q$ at $t_0$, such that for any $N \geq 0$, as $x \rightarrow + \infty$ we have
			\begin{equation}
				\int_{a} ^b	e^ {x p(t)}	q(t)dt	= e^ {x p(t_0)}	\left( \frac{c_0(t_0)}{\sqrt{x}}+	\frac{c_1(t_0)}{2! \, x^{3/2}}	+ \dots+	\frac{c_N(t_0)}{(2N)! \, x^{N+1/2}}+	\bigO{\frac{1}{x^{N+3/2}}}		\right)	\, .
			\end{equation}
	
		Moreover, $(c_N)_N$ can be computed as
			\begin{multline*}
			c_N(t_0)	=	\sqrt{\frac{2 \pi}{|	p''(t_0)|}}	\,	\sum_{k=0}^{2N}	\left( \begin{array}{c} 2N \\ k \end{array} \right) 	q^{(2N-k)}(t_0)	\\
		\sum_{m=0}^k	B_{k,m} \left(	\frac{p^{(3)}(t_0)	}{2.3}	, \dots		,\frac{p^{(k-m+3)}(t_0)		}{(k-m+2)(k-m+3)}		\right)	\frac{(2m+2N-1)!!}{| p''(t_0)|^{m+N}}.
			\end{multline*}
			where $B_{k,m}$ are the Bell polynomials defined above and $(2n+1)!!=1.3.5.\dots (2n+1)$.
	\end{theo}
	
\section*{Acknowlegments}
The authors wish to thank B. Bercu for valuable and fruitful discussions.
\bibliographystyle{plain}

\end{document}